\documentclass[11pt]{article}
\usepackage{comment}
\usepackage{amsfonts}
\usepackage{graphicx}
\usepackage{epsfig}
\usepackage{epstopdf}
\usepackage{multirow}
\usepackage{latexsym}
\usepackage{amssymb}
\usepackage{amsmath}
\usepackage{array}
\usepackage{xcolor}
\definecolor{cadmiumred}{rgb}{10,0,2}
\definecolor{red}{rgb}{1,0,0}
\definecolor{blue}{rgb}{0,0,1}
\definecolor{green}{rgb}{0,1,0}
\begin{document}
\newtheorem{t1}{Theorem}[section]
\newtheorem{d1}{Definition}[section]
\newtheorem{c1}{Corollary}[section]
\newtheorem{l1}{Lemma}[section]
\newtheorem{r1}{Remark}[section]
\newtheorem{ex}{Example}[section]
\newtheorem{re1}{Result}[section]
\newtheorem{co}{Counterexample}[section]
\newtheorem{p1}{Proposition}[section]
\title{Geometric and Harmonic Aging Intensity function and a Reliability Perspective}
\author{Subarna Bhattacharjee$^{1}$$^{}$\thanks{Corresponding author:~~
E-mail: subarna.bhatt@gmail.com, \\ $^{3}$ Sabana Anwar contributed in this article while she was a Project Assistant in Department of Mathematics, Ravenshaw University, Cuttack}
Ananda Sen$^{2}$, Sabana Anwar$^{3}$, Aninda K. Nanda$^{4}$\\
{\it $^{1,3}$ Department of Mathematics,Ravenshaw University, Cuttack-753003, Odisha, India}\\
{\it $^{2}$ Department of Family Medicine and Biostatistics, University of Michigan, Ann Arbor, MI, USA}\\
{\it $^{3}$ Institute of Mathematics and Applications,Bhubaneswar 751029, India}\\
{\it $^{4}$ Department of Statistics, Babasaheb Bhimrao Ambedkar University, Lucknow, UP, India}\\}
\date{June 13, 2024}
\maketitle
\begin{abstract}
In this paper, we introduce some new notions of aging based on geometric, harmonic means of failure rate and aging intensity function. We define a generalized version of aging functions called specific interval-average geometric hazard rate, specific interval-average harmonic hazard rate. We focus on some characterization results and their inter-relationships among the resulting non-parametric classes of distributions. Monotonic nature of so defined aging classes are exhibited by some well known probability distributions. Probabilistic orders based on these functions are taken up for further study. The work is illustrated through case studies and a simulated data having applications in reliability/survival analysis.
\end{abstract}
\section{Introduction}
In recent literature,  developing new aging functions and their analyses for its subsequent applications in various fields pertaining to study of aging phenomena is a thrust area among researchers. The present work is an attempt in this direction with emphasis on means of failure rate and other aging functions. To this end, here we place some of the crucial facts on means from vast literature.\\

In probabilistic framework of statitstics, the three types of mean, namely, arithmetic mean (AM), geometric mean (GM), and harmonic mean (HM) have been used extensively. There is a comprehensive review by Beebe (2023) (http://ftp.math.utah.edu/pub/tex/bib/agm.pdf) that detail their usage, mathematical properties, inter-relationship and extensions. Burk (1985) in the article entitled ``By all means" established the ordering $$HM<GM<LM<AM<RMS$$ where LM refers to the logarithmic mean (defined as $\frac{b-a}{\ln b-\ln a}$) and root mean square (RMS) (defined as $\sqrt{(a^2+b^2)/2}$. In mathematics and statistics, measures of central tendency offer a concise way to capture the overall characteristics of a dataset. While the pivotal position of central tendency in statistical infernce has been assumed by $SAM = \frac{1}{n} \sum_{i=1}^{n} x_i$, the sample version of geometric mean (SGM) defined as $SGM = (x_1x_2\ldots x_n)^{1/n}=\exp(\sum_{i=1}^n \frac{\ln x_i}{n})$ for $x>0$ has found numerous applications including environmental monitoring (e.g., acceptable level of contaminants in water quality and other immunologic information), infometrics, scientometrics (e.g., citation counts), finance (e.g., investment portfolio returns), nuclear medicine (e.g., tissue attenuation), ecology (e.g., growth rates in ecological population), groundwater hydrology, geoscience, geomechanics, machine learning (e.g., pattern recognition algorithm), chemical engineering (e.g., reaction rates), poverty and human development among others (c.f. Vogel R.M.(2020). A closely related measure is the sample version of the harmonic mean defined as $SHM = \frac{n}{\sum_{i=1}^{n}x_{i}^{-1}}$.\\

In nonparametric life-testing and reliability, the notion of aging has been a focal point of interest for several decades. While the failure rate $r(t)$ attempts to capture the aging behavior of a distribution, it acts as a poor comparator across distributions, especially if the rate is non-monotonic. Jiang (2003) introduced a quantitative tracking measure of aging, called aging intensity function (AI) defined by $L(t)=\frac{r(t)}{\frac{1}{t}\int_0^t r(u)du}$ for
$t>0.$ AI function written as

\begin{equation}
\label{new}L(t)=r(t)/A(t), \mbox{~where~} A(t)=\frac{1}{t}\int_0^t r(u)du,
\end{equation}
is the average failure rate. Because $L(t)$ expresses $r(t)$ in comparison to the average cumulative hazard at time $t$, it is a better tool to compare between distribution.\\

\hspace*{1cm}In this paper, we introduce some new functions which measure and explain the aging phenomenon of a system. The lifetime of any system (biological or non-living) having a well defined statistical distribution is represented by a random variable $X$. We give the following definition.\\

\hspace*{1cm}Let $X$ be a random variable having failure rate function $r_{X}(\cdot)$  such that $r_{X}(t)<\infty,$ for all $t>0.$ For $t>0,$ we define
\begin{equation}\label{new2}L^{G}_{X}(t)=\frac{r_{X}(t)}{G_{X}(t)}, ~L^{H}_{X}(t)=\frac{r_{X}(t)}{H_{X}(t)},\end{equation}
called the geometric aging intensity (GAI) and the harmonic aging intensity (HAI) of the random variable $X$, where $$G_{X}(t)=\exp\left(\frac{1}{t}\int_{0}^{t}\ln r_{X}(u)du\right), H_{X}(t)=\left(\frac{1}{t}\int_{0}^{t}\frac{1}{r_{X}(u)}du\right)^{-1}.$$ In view of the said functions, we
rename $L(t)$  as arithmetic aging intensity (also popularly known as aging intensity defined by Jiang et al. (2003)). Henceforth, unless otherwise required, we drop $X$ from all the expressions. \\
\hspace*{1cm}Readers may note that the functions $A(t), G(t)$ and $H(t)$ appearing in (\ref{new}) and (\ref{new2}) called arithmetic mean failure rate (AFR), geometric mean failure rate (GFR) and harmonic mean
failure rate (HFR) respectively were first coined by Roy and Mukherjee (1992). They defined the aging classes of lifetime distributions on the basis of monotonicity of GFR and
HFR. The very nomenclature of increasing (decreasing) geometric mean failure rate, increasing (decreasing) harmonic mean failure rate distributions were introduced by Roy and Mukherjee (1992) and they
denoted the resulting aging classes by IGFR (DGFR) and IHFR (DHFR) respectively. They proved that $IFR\subseteq IFRA\subseteq IGFR\subseteq IHFR$ and $DFR\subseteq DHFR\subseteq DGFR\subseteq DFRA$.
It may also be noted that Nanda et al. (2007), Bhattacharjee et al. (2013), Bhattacharjee et al. (2022), Giri et al. (2023), Szymkowiak (2018) and others worked on AI function and its properties.\\

\hspace*{1cm} The paper is organized as follows. Section 2 throws light on main results involving geometric and harmonic failure rates. Some characterization results and examples of parametric distributions are discussed for the aforementioned aging functions. Section 3  introduces new functions viz., specific interval-average geometric hazard rate, specific interval-average harmonic hazard rate and discusses their importance in reliability theory. In Section 4, we discuss some new probabilistic orders based on means of aging functions and obtain their inter-relationships. In Section 5, we implement the so obtained theoretical results on real life data and a simulated data. We note some significant observations with regard to bias and mean squared error of the estimators of aging functions. Section 6 gives the concluding remarks.

\section{Main Results}
In this section, we explore some characterizing properties of GAI and HAI functions. We group the characterization results into several categories starting with those based on nonparametric classification of distributions.
\subsection{Characterization of Aging Classes based on Geometric and Harmonic Aging Intensity function}
\hspace*{1cm} Lower and upper bounds of $L(t)$ for special classes of distributions have been established by Nanda et al (2007), Sunoj and Rasin (2018), Bhattacharjee et al.(2022).  The following theorem extends these results for $L^{G}(\cdot)$ and $L^{H}(\cdot)$.\\
\begin{t1}
For all $t>0,$
\label{impth}
\begin{enumerate}
\item[(i)] A random variable $X$ is IGFR (DGFR) if and only if $L^{G}(t)\geq (\leq)~1$.
\item[(ii)] A random variable $X$ is IHFR (DHFR) if and only if $L^{H}(t)\geq (\leq)~1$.
\item[(iii)]If $X$ is IFR (DFR) then $L^{G}(t)\geq( \leq)1$ and $L^{H}(t)\geq( \leq)1$.
\item[(iv)] If $X$ is IFRA then $L^{G}(t)> 1$  and $L^{H}(t)\geq 1$ for all $t>0.$
\item[(v)] If $X$ is DHFR then $L^{G}(t)\leq 1$ for all $t>0.$
\end{enumerate}
\end{t1}
{\bf Proof.}  A random variable $X$ is IGFR (DGFR) if and only  $G(t)$ is increasing (decreasing) in $t$. We know that for any function, say $g,$ $\exp(g(t))$ is increasing (decreasing) in $t$ if and only if $\frac{d}{dt}g(t)\geq (\leq)~0$ for all $t\geq 0$. This implies that $X$ is IGFR (DGFR) if and only if $\frac{d}{dt}\Big(\frac{1}{t}\int_{0}^{t} \ln r(u)du\Big)\geq(\leq)~0$ for all $t\geq 0,$ or equivalently $L^{G}(t)\geq(\leq)~1,$ thereby  proving $(i).$ Similarly, one can prove $(ii)$. To prove $(iii)$ we note that if $X$ is IFR then $L(t) \geq 1$ for all $t>0.$ Since, $L(t)\leq L^{G}(t)\leq L^{H}(t)$ for $t>0,$ it follows that $L^{G}(t)\geq 1$ and $L^{H}(t)\geq 1$ for all $t>0.$ $X$ is DFR is equivalent to the fact that $\ln r(u) \geq \ln r(t)$ for all $u\leq t.$ This gives $\int_{0}^{t}\ln r(u)du \geq t\ln r(t)$  i.e., $L^{G}(t)\leq 1$ for all $t> 0.$ Similarly, we prove that $L^{H}(t)\leq 1$ for all $t>0$ if $X$ is DFR. Parts (iv) and (v) follow from the aging class hierarchy provided in Roy and Mukherjee (1992). $\hfill\Box$\\

\hspace*{1cm}The next theorem gives lower and upper bounds of the functions $L, L^{G}$ and $L^{H}$ as stated below.
The following bounds are obtained for general random variables.
\begin{p1}
Let $X$ be a random variable with $0<r(t)<\infty,$ for $t\in R^{+}.$ Then
$$\frac{r(t)}{\sup_{t\in R^{+}} r(t)}\leq L(t)\leq L^{G}(t)\leq L^{H}(t)\leq  \frac{r(t)}{\inf_{t\in R^{+}} r(t)}.$$ Equality holds if and only if $X$ follows exponential distribution.
\end{p1}
\subsection{Characterization of parametric classes}
In this section we characterize some well-known distributions through the aging functions handled in this paper. \\

\hspace*{1cm}We know that the support of log-Weibull distribution is $(-\infty, \infty)$ and the corresponding survival function is $\bar{F}_{X}(t)=\exp\Big(-e^{\frac{t-a}{b}}\Big)$ with $a\in (-\infty, \infty), b>0.$ However, truncated log-Weibull distribution (cf. Giri et al. (2023) has positive real line as its support, i.e., $t\in (0,\infty)$ with survival function given by $$\bar{F}_{X}(t)=\frac{1}{e^{-e^{-\frac{a}{b}}}}\exp\Big(-e^{\frac{t-a}{b}}\Big).$$
The corresponding failure rate is
\begin{eqnarray}
r(t)&=&\frac{1}{b}e^{\frac{t-a}{b}}\nonumber\\
\label{rt}
&=&r(0)\Big(\frac{r(1)}{r(0)}\Big)^{t}
\end{eqnarray}
 In the next proposition, we characterize truncated log-Weibull distribution through $G(t).$
\begin{p1}
Let a non-negative random variable has failure rate $r(t)$ such that its first derivative exists. $G(t)=\sqrt{r(t)r(0)}$ for all $t\geq 0$ if and only if $X$ follows truncated log-Weibull distribution.
\end{p1}
{\bf Proof.} Clearly, $G(t)=\sqrt{r(t)r(0)}$ is equivalent to the fact that
\begin{displaymath}
e^{\frac{2}{t}\int_{0}^{t}\ln r(u)du}=r(t)r(0),
\end{displaymath}
giving \begin{equation}
\label{diff}
2 \int_{0}^{t}\ln r(u)du=t \ln r(t)+t \ln r(0).
\end{equation}
On differentiating (\ref{diff}), we obtain $\ln (\frac{r(t)}{r(0)})=t \frac{r^{'}(t)}{r(t)}.$ Thus, for $t\geq 0,$ \begin{equation}\label{diff1}\frac{d}{dt}r(t)=\frac{r(t)}{t}\ln (\frac{r(t)}{r(0)})\end{equation} Taking $p=\ln r(t),$ (\ref{diff1}) reduces to $$\frac{dp}{dt}-\frac{p}{t}=-\frac{1}{t}\ln r(0).$$ This gives $$r(t)=r(0)e^{kt},$$ where $k$ is an arbitrary constant. Clearly, $k=\ln \Big( \frac{r(1)}{r(0)}\Big).$ Thus, \begin{eqnarray}r(t)&=&r(0)e^{t \ln \Big( \frac{r(1)}{r(0)}\Big) }\nonumber\\
&=&r(0)\Big(\frac{r(1)}{r(0)}\Big)^{t}\nonumber,
\end{eqnarray}
which is failure rate of truncated Weibull distribution as given in (\ref{rt}). Conversely, if $X$ follows log-Weibull distribution then it is easy to verify that for all $t\geq 0,$
\begin{eqnarray}
G(t)&=&\frac{1}{b}\exp\Big(\frac{t-2a}{2b}\Big) \nonumber\\
&=& \frac{1}{b}\exp\Big(\frac{t-a}{2b}\Big)\exp\Big(-\frac{a}{2b}\Big)\nonumber\\
&=&\sqrt{r(t)r(0)},\nonumber
\end{eqnarray}
the last equality follows from (\ref{rt}). This completes the proof. $\hfill\Box$\\

\hspace*{1cm}Nanda et al. (2007) proved that for a random variable $X,$ $L(t)=1$ if and only if $r(t)$ is constant for all $t>0$. In other words, exponential distribution is
characterized by constant failure rate, and  AI function equals one.  They noted that $L(t)=\beta, t>0$ characterizes two-parameter Weibull distribution with shape parameter $\beta$.\\
\hspace*{1cm} The following proposition shows that GAI and HAI functions too characterize Weibull distribution.
\begin{p1}
Let $X$ be a random variable having GAI and HAI given by $L^{G}(\cdot)$ and $L^{H}(\cdot)$ respectively. For all $t>0,$
$L^{G}(t)=c$ if only if $X$ follows two-parameter Weibull distribution with shape parameter $\beta=\ln c+1$. Also, $L^{H}(t)=c$ if only if
    $X$ follows two-parameter Weibull distribution with shape parameter $\beta=2-c^{-1}$ such that $c>
    1/2.$ It may be noted that HAI characterizes only a particular family of Weibull distributions with shape parameter less than 2.$\hfill\Box$
\end{p1}
\subsection{Some illustrative examples}
\hspace*{1cm}Now, we look into  GAI and HAI of some well known parametric distributions. In each case, we mention the corresponding $AI$ function. Nanda et al. (2007) showed that monotonicity of FR is not transmitted to AI function. We observe a similar behavior for GAI and HAI functions.
\begin{ex}
\label{example1}
We assume $t\geq 0$ and all the parameters involved in the distributions are non-negative unless otherwise mentioned:
\begin{enumerate}
\item[(i)]Let $X$ be a random variable having Erlang distribution with failure rate $r(t)=\frac{\lambda^{2}t}{1+\lambda t}.$ $X$ is IFR. The corresponding $AI$ function is $L(t)=\frac{\lambda^{2}t^{2}}{(1+\lambda t)\big\{\lambda t-\ln (1+\lambda t)\big\}}$ and $X$ is DAI (cf. Nanda et al. (2007)).\\ Here, $L^{G}(t)=\big(1+\lambda t\big)^{1/\lambda t}$ and $L^{H}(t)=\frac{\lambda t}{(1+\lambda t)(\lambda t+\ln t)}.$  We note that, $\frac{d}{dt} L^{G}(t)=\frac{\big(1+\lambda t\big)^{1/\lambda t}}{\lambda t^{2}}\alpha(t)\leq 0$ where $\alpha(t)=\lambda t-(1+\lambda t)\ln (1+\lambda t)$ is decreasing in $t$ giving $\alpha(t)\leq \alpha(0)=0.$ This shows that $X$ is DGAI.\\
In the next few lines, we prove our following claim
$(a)$ if $t\leq e$ then $X$ is DHAI,
$(b)$ if $t\geq e$ and $\lambda \geq 1$ then $X$ is DHAI,
$(c)$ if  $\lambda \geq 1$ then $X$ is DHAI,
$(d)$ if  $\lambda \leq 1$ then $X$ has non-monotonic HAI.
We first take up the case of $t\leq e.$ Clearly, $\frac{1}{L^{H}(t)}=1+\lambda t+g(t)+\ln t,$ where $g(t)=\frac{\ln t}{\lambda t}.$ Here, $\frac{d}{dt}g(t)=\frac{\lambda (1-\ln t)}{\lambda^{2}t^{2}}\geq 0$  if $t\leq e.$ Thus, $\frac{1}{L^{H}(t)}$ is non-decreasing in $t$ for $t\leq e.$ In other words, $X$ is DHAI for $t\leq e,$ which proves $(a)$.\\ We take up the case when $t\geq e$ and $\lambda \geq 1.$ We note that, \begin{equation}\label{firsteq}\frac{d}{dt}\frac{1}{L^{H}(t)}=\frac{\lambda^{2}t^{2}+\lambda t+1-\ln t}{\lambda t^{2}}=\frac{g_2(t)}{\lambda t^{2}}, (say),
\end{equation}
 and $\frac{d}{dt}\Big(\lambda^{2}t^{2}+\lambda t+1-\ln t\Big)=\frac{g_3(t)}{t},$ where $g_2(e)=\lambda^{2}e^{2}+\lambda e\geq 0$ and $g_3(t)=2\lambda^{2}t^{2}+\lambda t-1, (say).$ Here, $g_3(t) \geq e\lambda(1+2e\lambda)-1\geq0$ if $\lambda \geq 1.$
Thus, $g_2(t)=\lambda^{2}t^{2}+\lambda t+1-\ln t \geq g_2(e)=\lambda^{2}e^{2}+\lambda e\geq 0$ for $t\geq e$ and $\lambda \geq 1$ implies $X$ is DHAI for $t\geq e$ and $\lambda \geq 1.$ This proves $(b)$. Clearly, $(c)$ follows from $(a)$ and $(b).$\\
We now take up case $(d)$. If $\lambda<1$ and $t\geq e$, $L^{H} (t)$ needs further investigation. In (\ref{firsteq}), note that
\begin{eqnarray}
g_2(t)&=&\lambda^{2}t^{2}+\lambda t+1-\ln t
=(\lambda t+\frac{1}{2})^{2}-\ln t+\frac{3}{4}
\leq 0,\nonumber
\end{eqnarray}
if $\lambda\leq -\frac{1}{2t}+\frac{1}{2t}\sqrt{4\ln t-3}.$ Since, $t\geq e$, it follows that $g_2(t)\leq 0$ if $\lambda\leq \frac{\sqrt{4\ln t-3}-1}{2e}.$ Further, $\lambda\leq 1$ gives $t\leq e^{1+e+e^{2}}$. Thus, for $t\in (e,e^{1+e+e^{2}})$ and $\lambda<1$, $X$ is IHAI. We conclude that $L^{H}(t)$ is non-monotonic if $\lambda \leq1,$ thereby proving $(d).$
\item[(ii)] Let $X$ be a random variable having Uniform distribution with support $[a,b],$ and failure rate $r(t)=1/(b-t)$ for $a<t<b.$ Clearly, $X$ is IFR. We now study the monotonic property of AI, GAI and HAI functions. Here, $A(t)=\frac{1}{t-a}\Big(\int_{a}^{t}r(u)du\Big)=\frac{1}{t-a}\ln \Big(\frac{b-a}{b-t}\Big)$ and consequently, for $a<t<b,$ $$L(t)=\frac{t-a}{(b-t)\big(\ln (b-a)-\ln (b-t)\big)}.$$
    Note that,
    $\frac{d}{dt}L(t)=\frac{(b-t)\phi(t)}{\Big\{(b-t)\big(\ln (b-a)-\ln (b-t)\big)\Big\}^{2}},$
    where $\phi(t)=\ln (\frac{b-a}{b-t})-\frac{t-a}{b-t}.$ Since, $\frac{d}{dt}\phi(t)=\frac{a-t}{(b-t)^{2}}\leq 0,$ it follows that $\phi(t)$ is decreasing in t for $a<t<b.$ This gives, $\phi(t)\leq \phi(a),$ i.e., $\phi(t)\leq 0$ for $a<t<b.$ Thus, $\frac{d}{dt}L(t)\leq 0$ proving $X$ is DAI.
   We get $G(t)=\frac{e}{(b-t)^{\frac{t-b}{t-a}}(b-a)^{\frac{b-a}{t-a}}}$ and $L^{G}(t)=\frac{1}{e}\Big(\frac{b-t}{b-a}\Big)^{\frac{a-b}{t-a}}$ for $a<t<b.$ Also, $\frac{d}{dt}L^{G}(t)=\frac{1}{(a-t)^{2}}\Big(\frac{t-b}{a-b}\Big)^{\frac{t+b-2a}{a-t}}\alpha(t),$ where $\alpha(t)=\Big\{t-a+(b-t)\ln(\frac{t-b}{a-b})\Big\}.$ Since, $\frac{d}{dt}\alpha(t)=-\ln (\frac{t-b}{a-b})\geq 0,$ we get $\alpha(t)\geq 0$ and we arrive at the conclusion that $X$ is IGAI.\\
    Similarly, we obtain, $L^{H}(t)=1+\frac{t-a}{2(b-t)}$ and observe that $X$ is IHAI.
\item[(iii)]For Rayleigh distribution with failure rate $r(t)=a+bt,$ we have $L^{G}(t)=\frac{e}{(1+\frac{b}{a}t)^{a/bt}}$. $L^{H}(t)=\big(1+\frac{a}{bt}\big)\ln\big(1+\frac{bt}{a}\big).$ As noted by Nanda et al. (2007), $L(t)=\frac{a+bt}{a+bt/2}$ and $X$ is IAI. We note that $X$ is IGAI and IHAI.
\item[(iv)]For Pareto distribution with failure rate $r(t)=\frac{a}{t},$ $t\geq k$, we have $L^{G}(t)=\frac{1}{et}\Big(\frac{t^{t}}{k^{k}}\Big)^{\frac{1}{t-k}}$ and $L^{H}(t)=\frac{t+k}{2t}$. However, Nanda et al. (2007) showed that $L(t)=\frac{1}{\ln t-\ln k}$ and $X$ is DAI. We note that, $X$ is DGAI and DHAI.
\end{enumerate}
\end{ex}

\subsection{Some system properties in terms of GAI and HAI}
\hspace*{1cm}Bhattacharjee et al. (2013) showed that if a series system is formed by $n$ independent components with lifetimes denoted by $X_{i}$ for $1\leq i \leq n$ and
system lifetime $X,$ then its corresponding AI function $L_{X}(t)$ is bounded between $\min_{1\leq i \leq n} L_{X_{i}}(t)$  and $\max_{1\leq i \leq n} L_{X_{i}}(t)$, i.e.,
$\min_{1\leq i \leq n} L_{X_{i}}(t)\leq L_{X}(t)\leq \max_{1\leq i \leq n} L_{X_{i}}(t)$ for all $t>0.$ Here, we prove a similar result for $L^{G}(\cdot).$\\
\hspace*{1cm}An extension of a result due to Hardy et al. (2020) is given in the following Lemma and is used in upcoming theorem. We give an outline of the proof to relish its mathematical  rigour.
\begin{l1}\label{little}
$G(\sum_{1\leq i\leq n}f_{i})\geq \sum_{1\leq i\leq n}G(f_{i}),$ where $f_i>0$ for all $i=1,2,\ldots,n.$
\end{l1}
{\bf Proof.} Note that, \begin{eqnarray}
\frac{G(f_{i})}{G(\sum_{1\leq i\leq n}f_{i})}&=&\frac{\exp \Big(\frac{1}{t}\int_{0}^{t}\ln f_{i} du)}{\exp(\frac{1}{t}\int_{0}^{t}\ln (\sum_{1\leq i\leq n}f_{i}) du\Big)}\nonumber\\
&=&\exp \Big(\frac{1}{t}\int_{0}^{t}\ln \big(\frac{f_{i}}{\sum_{1\leq i\leq n}f_{i}}\big)du\Big)=
G\Big(\frac{f_{i}}{\sum_{1\leq i\leq n}f_{i}}\Big)
\leq A\Big(\frac{f_{i}}{\sum_{1\leq i\leq n}f_{i}}\Big),\nonumber
\end{eqnarray}
giving $\frac{\sum_{1\leq i\leq n}G(f_{i})}{G(\sum_{1\leq i\leq n}f_{i})}\leq \sum_{1\leq i\leq n}A(\frac{f_{i}}{\sum_{1\leq i\leq n}f_{i}})=1$. This proves that $G(\sum_{1\leq i\leq n}f_{i})\geq \sum_{1\leq i\leq n}G(f_{i})$.$\hfill\Box$ \\

\hspace*{1cm}Now, we mention a mathematical tool from Bhattacharjee et al. (2013), which is to be used in following theorem.
\begin{l1}
\label{bhatta}
Let $p_i > 0, q_i > 0,$ for $i = 1, 2,\ldots, k.$ Then
$$\min_{1\leq i\leq k}\big(\frac{p_i}{q_i}
\big) \leq \frac{\big(\sum_{1\leq i\leq k}p_i\big)}{\big(\sum_{1\leq i\leq k}q_i\big)}\leq \max_{1\leq i\leq k}\big(\frac{p_i}{q_i}
\big) $$
\end{l1}

\hspace*{1cm}The next theorem gives a upper bound of GAI function, $L^{G},$ of a series system formed by $n$ independent components.
\begin{t1}
  If a series system is formed by $n$ independent components with lifetimes denoted by $X_{i}$ for $1\leq i \leq n$ and system lifetime $X,$ then its corresponding geometric AI function,
  $L^{G}_{X}(t)$ is bounded above by $\max_{1\leq i \leq n} L^{G}_{X_{i}}(t)$, i.e., $L^{G}(t)\leq \max_{1\leq i \leq n} L^{G}_{X_{i}}(t)$ for all $t>0.$
\end{t1}
{\bf Proof.} Here, $r_{X}(t)=\sum_{i=1}^{n}r_{X_{i}}(t),$ so
\begin{eqnarray}
L^{G}(t)&=& \frac{\sum_{1\leq i\leq n}r_{X_{i}}(t)}{\exp\Big\{\frac{1}{t}\int_{0}^{t}\ln \Big(\sum_{1\leq i\leq n}r_{X_{i}}(t)\Big) du\Big\}}\nonumber\\
&=&\frac{\sum_{1\leq i\leq n}r_{X_{i}}(t)}{G(\sum_{1\leq i\leq n}r_{X_{i}}(t))}
\leq\frac{\sum_{1\leq i\leq n}r_{X_{i}}(t)}{\sum_{1\leq i\leq n}G(r_{X_{i}}(t))}
\leq \max_{1\leq i\leq n} L^{G}_{X_{i}}(t).\nonumber
\end{eqnarray}
The last two inequalities follow due to Lemma \ref{little} and upper bound given in Lemma \ref{bhatta}. This completes the proof. $\hfill\Box$\\

\hspace*{1cm} We know that AFR of a series system formed by finite number of independent components is the sum of the AFR of individual components, i.e., if $X=\min_{1\leq i\leq n}X_{i}$ and $X_{i}'s$ are independent random variables then $A_{X}(t)=\sum_{i=1}^{n}A_{X_{i}}(t)=n~DAM \{A_{X_{i}}(t):1\leq i\leq n\}$ where $DAM(c_{i}:1\leq i\leq n)$ denotes discrete-arithmetic mean of  $n$ non-negative quantities $c_{i}'s.$ This motivates one to explore about GFR and HFR of a series system comprised of independent components. The following theorem is interesting because of its mathematical coherence. This result gives a lower bound of $G_{X}(\cdot)$ and $H_{X}(\cdot)$ of a series system formed by $n$ independent components.
\begin{t1}
\label{series}
Let us consider a series system with lifetime denoted by $X$ which is formed by $n$ independent components with lifetimes $X_{i}$ for $1\leq i \leq n$, i.e., $X=\min_{1\leq i\leq n}X_{i}$. Then, for all $t>0,$
  $$ G_{X}(t)\geq n \Big(\prod_{i=1}^{n}G_{X_{i}}(t)\Big)^{\frac{1}{n}}=n~DGM \{G_{X_{i}}(t):1\leq i\leq n\}$$ and $$H_{X}(t)\geq n^{2} \Big(\sum_{i=1}^{n}\frac{1}{H_{X_{i}}(t)}\Big)^{-1}=n~DHM \{H_{X_{i}}(t):1\leq i\leq n\},$$ where $$DGM \{a_{i}:1\leq i \leq n, a_i >0\}=\Big(\prod_{i=1}^{n}a_{i}\Big)^{1/n}$$  and $$DHM \{b_{i}:1\leq i\leq n, b_i >0\}=n\Big(\sum_{i=1}^{n}\frac{1}{b_{i}}\Big)^{-1}$$ represent (discrete) geometric mean and (discrete) harmonic mean of $n$ number of finite non-negative numbers $a_{i}'s$ and $b_{i}'s$ respectively.
\end{t1}
{\bf Proof.} Since, $\sum_{i=1}^{n}r_{X_{i}}(t)\geq n\Big(\prod_{i=1}^{n}r_{X_{i}}(t)\Big)^{1/n},$ we get $$\frac{1}{t}\int_{0}^{t}\ln \Big(\sum_{i=1}^{n}r_{X_{i}}(u)\Big)du\geq \ln n +\frac{1}{n}\sum_{i=1}^{n}\frac{1}{t}\Big(\int_{0}^{t}\ln r_{X_{i}}(u)du\Big).$$ This implies, \begin{eqnarray}G_{X}(t)=\exp\Big\{\frac{1}{t}\int_{0}^{t}\ln \Big(\sum_{1\leq i \leq n}r_{X_{i}}(u)\Big)du\Big\}&\geq& n \exp\Big\{\frac{1}{n}\sum_{i=1}^{n}\frac{1}{t}\Big(\int_{0}^{t}\ln r_{X_{i}}(u)du\Big)\Big\}\nonumber\\
&=& n \exp\Big\{\frac{1}{n}\sum_{i=1}^{n}\ln G_{X_{i}}(t)\Big\}\nonumber\\
&=& n\Big(\prod_{i=1}^{n}G_{X_{i}}(t)\Big)^{\frac{1}{n}}
=n~DGM \{G_{X_{i}}(t):1\leq i\leq n\}.\nonumber
\end{eqnarray}
To obtain a bound for HFR of $X$, we recall that $\sum_{i=1}^{n}r_{X_{i}}(t)\geq n^{2}\Big(\sum_{i=1}^{n}\frac{1}{r_{X_{i}}(t)}\Big)^{-1}.$
Here,
\begin{eqnarray}
\Big(H_{X}(t)\Big)^{-1}=\frac{1}{t}\int_{0}^{t}\frac{1}{r_{X}(u)}du
&=&\frac{1}{t}\int_{0}^{t}\frac{1}{\sum_{i=1}^{n}r_{X_{i}}(u)}du\nonumber\\
&\leq&\frac{1}{n^{2}}\Big\{\frac{1}{t}\int_{0}^{t}\sum_{i=1}^{n}\frac{1}{r_{X_{i}}(u)}du\Big\}\nonumber\\
&=&\frac{1}{n^{2}}\sum_{i=1}^{n}\frac{1}{t}\int_{0}^{t}\frac{1}{r_{X_{i}}(u)}du
=\frac{1}{n^{2}}\sum_{i=1}^{n}\frac{1}{H_{X_{i}}(t)}\nonumber
\end{eqnarray}
giving $$H_{X}(t)\geq n^{2} \Big(\sum_{i=1}^{n}\frac{1}{H_{X_{i}}(t)}\Big)^{-1}=n~DHM \{H_{X_{i}}(t):1\leq i\leq n\}.$$ This completes the proof.$\hfill\Box$\\

\hspace*{1cm}A sharper lower and upper bound of AFR, GFR and HFR of a series system constituted by $n$ independent components is given in following remark. The symbols used have their usual meaning as described in this paper.
\begin{r1}
Let us consider a series system with lifetime denoted by $X$ which is formed by $n$ independent components with lifetimes $X_{i}$ for $1\leq i \leq n$, i.e., $X=\min_{1\leq i\leq n}X_{i}$. Then, for all $t>0,$
$$nDAM\{A_{X_{i}}(t):1\leq i\leq n\}=A_{X}(t)\geq G_{X}(t)\geq H_{X}(t)\geq n~DHM \{H_{X_{i}}(t): 1\leq i\leq n\}.$$
\end{r1}
\section{Some new notions of aging: specific interval-average geometric hazard rate, specific interval-average harmonic hazard rate}
\hspace*{1cm}Before we  begin this section, readers may note that the notations used here, namely, $A(t,s),$ $G(t,s)$ and $H(t,s)$ shall be separately dealt with and are not related to the notations that are referred in earlier sections.\\
\hspace*{1cm}Bryson and Siddiqui (1969) introduced some notions of aging, called specific aging factor, denoted by $A(t,s)$ and the specific interval-average hazard rate, denoted by $H(t,s)$ of a system at time $t$ specific with respect to a positive time parameter $s$. They defined $A(t,s)=\frac{\overline{F}(t)\overline{F}(s)}{\overline{F}(t+s)},$  and $H(t,s)=\frac{1}{t}\int_{s}^{s+t}r(u)du$ for all $t,s\geq 0.$ \\\hspace*{1cm}Inquisitive readers may study the applications of $A(t,s)$ in comparison of two systems with different chronological age, older system having age $t$ and other one is new with chronological age zero but both having the same survival function, say $\bar{F}(\cdot)$. Clearly, $A(t,s)$ is the ratio of survival probabilities of new and older systems, i.e., it is the ratio of $\overline{F}(s)$ (i.e., the survival probability that the new system will survive for at least a duration of $s$ units ) and $\frac{\overline{F}(t+s)}{\overline{F}(t)}$ (i.e., the survival probability that older system will survive for same $s$ duration, given its prior survival up to time $t$).  Both these quantities are defined with respect to a positive time parameter $s$. Clearly, $H(t,s)$ is a generalization of hazard rate average $H(t,0)=\frac{1}{t}\int_{0}^{t}r(u)du$ which was introduced by Birnbaum et al. (1966). Bryson and Siddiqui (1969) proved that $H(t,s)$ is increasing in $t$ if and only if $r(t)$ is increasing in $t$. \\
\hspace*{1cm}In this section, we introduce the generalized version of GFR and HFR termed as specific interval-average geometric hazard rate and specific interval-average harmonic hazard rate respectively as given in the following definition. The importance of these generalized functions can be seen in upcoming Theorems \ref{gm} and \ref{gmh}.
\begin{d1}
The specific interval-average geometric hazard rate of a non-negative random variable $X$ having hazard rate $r(\cdot)$ is $$GM_{X}(t,s)=\exp\left(\frac{1}{t}\int_{s}^{s+t}\ln r_{X}(u)du\right), \mbox{~ for~} t,s\geq0.$$
\end{d1}
\begin{d1}
The specific interval-average harmonic hazard rate of a non-negative random variable $X$ having hazard rate $r(\cdot)$ is $$HM_{X}(t,s)=\left(\frac{1}{t}\int_{s}^{t+s}\frac{1}{r_{X}(u)}du\right)^{-1}, \mbox{~ for~} t,s\geq 0.$$
\end{d1}
We drop $X$ from $GM_{X}(t,s)$ and $HM_{X}(t,s)$ and simply write $GM(t,s)$ and $HM(t,s)$ respectively.\\

\hspace*{1cm} Readers would like to immediately look into the significance of above functions. To this end, we focus on residual lifetime of a random variable $X$ denoted by $R^{X}_{x}=(X-x\mid X>x)$ for $x>0$. We note that the AFR of residual life-time $R^{X}_{x}$ is its specific aging factor, $i.e.,$  $A_{R^{X}_{x}}(t)=A(t,x),$ GFR of residual lifetime $R^{X}_{x}$ is equal to its specific interval-average geometric hazard rate, i.e., $G_{R^{X}_{x}}(t)=GM(t,x),$ and HFR of residual life-time $R^{X}_{x}$ is equal to its specific interval-average harmonic hazard rate, i.e.,  $H_{R^{X}_{x}}(t)=HM(t,x).$ As a result, it follows that GAI and HAI of residual lifetime $R^{X}_{x}$ are $L^{G}_{R^{X}_{x}}(t)=\frac{r(x+t)}{GM(t,x)}$ and $L^{H}_{R^{X}_{x}}(t)=\frac{r(x+t)}{HM(t,x)}$ for all $t>0.$\\

\hspace*{1cm} In upcoming theorem, we establish an equivalent condition of monotonic increasing hazard rate in terms of specific interval-average geometric hazard rate.
\begin{t1}
\label{brynew1}
Let $h(t)$ be integrable, with no more than finitely many discontinuities in any finite interval. Then $h(t)$ is monotone increasing for all $t>0$ if and only if \begin{equation}\label{eq1}h(s)\leq GM(t,s)\leq h(s+t)\end{equation}for all $s\geq 0, t\geq 0.$
\end{t1}
{\bf Proof.} Under the hypothesis, clearly $h(t)$ is monotone increasing for all $t>0$ if and only if \begin{equation}\label{brysi}h(s)\leq \frac{1}{t}\int_{s}^{s+t} h(x)dx\leq h(s+t)\end{equation} for all $s\geq 0, t\geq 0,$ (cf. Bryson and Siddiqui (1969)). Since $h(t)$ is monotonic increasing in $t$ is equivalent to $\ln h(t)$ being monotonic increasing, we replace $h(t)$ by $\ln h(t)$ in (\ref{brysi}). Thereby, we conclude that (\ref{brysi}) is equivalent to   $$\ln h(s)\leq \frac{1}{t}\int_{s}^{s+t} \ln h(x)dx\leq \ln h(s+t)$$ for all $s\geq 0, t\geq 0.$ This completes the proof.$\hfill\Box$\\

\hspace*{1cm} To have the counterpart of HFR, we state the following theorem.
\begin{t1}
\label{brynewhq}
Let $h(t)$ be integrable, with no more than finitely many discontinuities in any finite interval. Then $h(t)$ is monotone increasing for all $t>0$ if and only if \begin{equation}\label{eq2}h(s)\leq HM(t,s)\leq h(s+t)\nonumber\end{equation}for all $s\geq 0, t\geq 0.$
\end{t1}
{\bf Proof.} It is easy to prove that $h(t)$ is monotone increasing for all $t>0$ if and only if $$ h(s)\leq \Big(\frac{1}{t}\int_{s}^{s+t} \frac{1}{h(x)}dx\Big)^{-1}\leq h(s+t)$$ for all $s\geq 0, t\geq 0.$ This proves the theorem.$\hfill\Box$\\

\hspace*{1cm}A motivation of the above defined aging functions, namely  specific interval-average geometric hazard rate and specific interval-average harmonic hazard rate is highlighted in the upcoming two theorems.
\begin{t1}
\label{gm}
A random variable $X$ has increasing (decreasing) hazard rate $r(\cdot)$ if and only if $$GM(t_2,s)\geq(\leq)GM(t_1,s)$$ for all $t_2\geq t_1\geq 0,$ and $s\geq 0.$
\end{t1}
{\bf Proof.} For $t_2\geq t_1\geq 0,$
\begin{eqnarray}
&&GM(t_2,s)-GM(t_1,s)\nonumber\\&=&\exp\Big(\frac{1}{t_2}\int_{s}^{s+t_2}\ln r_{X}(u)du\Big)-GM(t_1,s)\nonumber\\
&=&\tiny{\exp\Big(\frac{1}{t_2}\int_{s}^{s+t_1}\ln r_{X}(u)du\Big)\exp\Big(\frac{1}{t_2}\int_{s+t_1}^{s+t_2}\ln r_{X}(u)du\Big)-GM(t_1,s)}\nonumber\\
&=&\exp\Big(\frac{1}{t_2}\int_{s}^{s+t_1}\ln r_{X}(u)du\Big)\exp\Big(\frac{1}{t_2-t_1}\big(\frac{t_2-t_1}{t_2}\Big)\int_{s+t_1}^{s+t_2}\ln r_{X}(u)du\Big)-GM(t_1,s)\nonumber\\
\label{eq31}&=&\Big((GM(t_1,s))^{\frac{t_1}{t_2}} ~GM(t_2-t_1,s+t_1)\Big)^{1-\frac{t_1}{t_2}}-GM(t_1,s)
\end{eqnarray}
If $h(t)$ is increasing in $t$ then using (\ref{eq1}) in (\ref{eq31}), we find that $GM(t_2-t_1,s+t_1))^{1-\frac{t_1}{t_2}}\geq h(s+t_1)^{1-\frac{t_1}{t_2}}.$ As a result, we obtain \begin{eqnarray}GM(t_2,s)-GM(t_1,s)&\geq& (GM(t_1,s))^{\frac{t_1}{t_2}} h(s+t_1)^{1-\frac{t_1}{t_2}}-GM(t_1,s))\nonumber\\
\label{eq3}&\geq& GM(t_1,s)\Big[\Big(\frac{h(s+t_1)}{GM(t_1,s)}\Big)^{1-\frac{t_1}{t_2}}-1\Big]
\end{eqnarray}
We observe that $\Big(\frac{h(s+t_1)}{GM(t_1,s)}\Big)^{1-\frac{t_1}{t_2}}\geq 1$ as $1-\frac{t_1}{t_2}\geq 0$ and $\frac{h(s+t_1)}{GM(t_1,s)}\geq 1.$ Hence from (\ref{eq3}), it follows that $GM(t_2,s)-GM(t_1,s)\geq 0$ for all $t_2\geq t_1\geq 0,$ and for all $s\geq 0$ if $h(t)$ is increasing in $t$. One can prove the converse part in a similar manner. This proves the theorem.$\hfill\Box$
\begin{t1}
\label{gmh}
A random variable has increasing (decreasing) hazard rate if and only if $$HM(t_2,s)\geq(\leq)HM(t_1,s)$$ for all $t_2\geq t_1\geq 0,$ $s\geq 0.$
\end{t1}
{\bf Proof.} Let $h(x)$ be monotone increasing and $t_2\geq t_1\geq 0.$ Then
\begin{eqnarray}
\frac{1}{HM(t_2,s)}-\frac{1}{HM(t_1,s)}&=&\frac{1}{t_2}\int_{s}^{s+t_2}\frac{1}{h(u)}du-\frac{1}{HM(t_1,s)}\nonumber\\
&=&\frac{1}{t_2}\Big(\int_{s}^{s+t_1}\frac{1}{h(u)}du+\int_{s+t_1}^{s+t_2}\frac{1}{h(u)}du\Big)-\frac{1}{HM(t_1,s)}\nonumber\\
&=&\Big(\frac{t_1}{t_2}\Big)\frac{1}{HM(t_1,s)}+\frac{t_2-t_1}{t_2}\frac{1}{HM(t_2-t_1,s+t_1)}-\frac{1}{HM(t_1,s)}\nonumber\\
&=&\big(\frac{t_1-t_2}{t_2}\big)\frac{1}{HM(t_1,s)}+\big(\frac{t_2-t_1}{t_2}\big)\frac{1}{HM(t_2-t_1,s+t_1)}\nonumber\\
&\leq&\big(\frac{t_1-t_2}{t_2}\big)\frac{1}{h(s)}+\big(\frac{t_2-t_1}{t_2}\big)\frac{1}{h(s+t_1)}
=\big(\frac{t_1-t_2}{t_2}\big)\frac{h(s+t_1)-h(s)}{h(s)h(s+t_1)}
\leq 0,\nonumber
\end{eqnarray}
proving one part of the theorem. We can prove the converse part easily. This completes the proof.$\hfill\Box$\\

\hspace*{1cm}Theorem \ref{gm} and Theorem \ref{gmh} are restated in the following corollary.
\begin{c1}
$X$ is IFR (DFR) if and only if $GM(t,s)$ is increasing (decreasing) in $t$ for all $s\geq 0.$ Similarly, $X$ is IFR (DFR) if and only if $HM(t,s)$ is increasing (decreasing) in $t$ for all $s\geq 0.$$\hfill\Box$
\end{c1}

\hspace*{1cm}The aforementioned two theorems help us to immediately infer that the conditions $r(t), GM(t,s)$ and $HM(t,s)$  increasing (decreasing) in $t$ for all $s\geq 0,$ are equivalent. However, the importance of increasing (decreasing) aging classes based on $GM(t,s)$ and $HM(t,s)$ is due to the fact that they are generalization of  geometric failure and harmonic failure rates respectively.
\section{Probabilistic Order based on means of aging functions}
\hspace*{1cm}The role of stochastic orders are widely accepted in various fields by scientists.
In this section, we give a brief account of some new orders based on aging functions defined in the present work.\\\hspace*{1cm} Let $X$ and $Y$ be two random variables with failure rates FR, arithmetic mean failure rates AFR, geometric failure rates GFR, harmonic failure rates HFR, geometric aging intensity GAI, harmonic aging intensity HAI, given by $r_{X}(t), r_{Y}(t);$ $A_{X}(t), A_{Y}(t);$ $G_{X}(t), G_{Y}(t);$ $H_{X}(t), H_{Y}(t);$ $L^{G}_{X}(t), L^{G}_{Y}(t);$$L^{H}_{X}(t), L^{H}_{Y}(t);$ respectively. The following definition gives a detailed account of the newly introduced orders namely $AFR,GFR,HFR,GAI$ and $HAI$.
\begin{d1}
A random variable $X$ is said to be smaller than another random variable $Y$ in
\begin{enumerate}
\item[(i)] arithmetic mean failure rate (denoted by $X\leq_{AFR} Y$) if $A_{X}(t)\geq A_{Y}(t);$ for all $t>0,$
\item[(ii)] geometric mean failure rate (denoted by $X\leq_{GFR} Y$) if $G_{X}(t)\geq G_{Y}(t);$ for all $t>0,$
\item[(iii)] harmonic mean failure rate (denoted by $X\leq_{HFR} Y$) if $H_{X}(t)\geq H_{Y}(t);$ for all $t>0,$
\item[(iv)]  geometric aging intensity (denoted by $X\leq_{GAI} Y$) if $L^{G}_{X}(t)\geq L^{G}_{Y}(t);$ for all $t>0,$
\item[(v)] harmonic aging intensity (denoted by $X\leq_{HAI} Y$) if $L^{H}_{X}(t)\geq L^{H}_{Y}(t);$ for all $t>0.$$\hfill\Box$
\end{enumerate}
\end{d1}

The reflexive, commutative and antisymmetric properties of GAI and HAI order are given below.
\begin{t1}
\begin{enumerate}
\item[(i)] $X\leq_{GAI}X$  and $X\leq_{HAI}X.$
\item[(ii)] If $X\leq_{hr}Y$ then $X\leq_{AFR}Y.$ If $X\leq_{hr}Y$ then $X\leq_{HFR}Y.$
\item[(iii)] If $X\leq_{GAI}Y$ and $Y\leq_{GAI}Z$ then $X\leq_{GAI}Z.$ If $X\leq_{HAI}Y$ and $Y\leq_{HAI}Z$ then $X\leq_{HAI}Z.$
\item[(iv)]  If $X\leq_{GAI}Y$ and $Y\leq_{GAI}X$ then $X$ and $Y$ have proportional failure rates. If $X\leq_{HAI}Y$ and $Y\leq_{HAI}X$ then $X$ and $Y$ have proportional failure rates.
\end{enumerate}
\end{t1}
{\bf Proof.} Proofs of (i), (ii) and (iii) are straightforward. If $X\leq_{GAI}Y$ and $Y\leq_{GAI}X$ then $L^{G}_{X}(t)=L^{G}_{Y}(t),$  which on simplification gives $t\frac{r_{Y}(t)}{r_{X}(t)}\frac{d}{dt}\Big(\frac{r_{Y}(t)}{r_{X}(t)}\Big)=0$ for all $t>0$. This implies $r_{Y}(t)=cr_{X}(t)$ for all $t>0.$ If $X\leq_{HAI}Y$ and $Y\leq_{HAI}X$ then $L^{H}_{X}(t)=L^{H}_{Y}(t),$ gives \begin{equation}\label{hr1}r_{X}(t)\int_{0}^{t}\frac{1}{r_{X}(u)}du=r_{Y}(t)\int_{0}^{t}\frac{1}{r_{Y}(u)}du.\end{equation} On differentiating, we get
\begin{equation}
\label{hr2}
\Big(\frac{d}{dt}r_{X}(t)\Big)\Big(\int_{0}^{t}\frac{1}{r_{X}(u)}du\Big)=\Big(\frac{d}{dt}r_{Y}(t)\Big)\Big(\int_{0}^{t}\frac{1}{r_{Y}(u)}du\Big).
\end{equation}
From (\ref{hr1}) and (\ref{hr2}), we find that $\frac{r_{X}(t)}{r_{Y}(t)}=\frac{r_{X}^{'}(t)}{r^{'}_{Y}(t)}$ where $'$ represents differentiation with respect to $t.$  This implies, $\frac{d}{dt}\Big(\frac{r_{Y}(t)}{r_{X}(t)}\Big)=0.$ This proves $(iv)$. $\hfill\Box$\\

\hspace*{1cm}It is worthwhile to note from Sengupta and Deshpande (1994) and Rowell and Siegrist (1998) that monotonicity of the ratio $r_{X}(t)/r_{Y}(t)$ for all $t\geq 0$ is an equivalent condition to say that a random variable $X$ is aging faster than $Y$, denoted by $X\leq_{AF}Y.$
Nanda et al. (2007) noted that  $X\leq_{AF}Y$ then $X\leq_{AI}Y$.\\

\hspace*{1cm}The next theorem gives an equivalent condition of $GAI$ order. We also prove that $AF$ order is stronger than $GAI$ order and give a sufficient condition for $GAI$ order.
\begin{t1}
For two random variables $X$ and $Y,$
\begin{enumerate}
\item[(i)] $X\leq_{GAI}Y$ if and only if $\ln \Big(\frac{r_{X}(t)}{r_{Y}(t)}\Big)\geq \frac{1}{t}\int_{0}^{t}\ln \Big(\frac{r_{X}(u)}{r_{Y}(u)}\Big)du$, for $t>0$
\item[(ii)] If $\ln \Big(\frac{r_{X}(t)}{r_{Y}(t)}\Big)$ is increasing in $t>0$ then $X\leq_{GAI}Y.$
\item[(iii)]If $X\leq_{AF}Y$ then $X\leq_{GAI}Y$.
\end{enumerate}
\end{t1}
{\bf Proof.} We note that $X\leq_{GAI}Y$ if and only if $r_{X}(t)\exp\big(\frac{1}{t}\int_{0}^{t}\ln r_{Y}(u)du\big) \geq r_{Y}(t)\exp\big(\frac{1}{t}\int_{0}^{t}\ln r_{X}(u)du\big)$. Thus, $X\leq_{GAI}Y$ is equivalent to the fact that $\ln \Big(\frac{r_{X}(t)}{r_{Y}(t)}\Big)\geq \frac{1}{t}\int_{0}^{t}\ln \Big(\frac{r_{X}(u)}{r_{Y}(u)}\Big)du$ for all $t>0.$ This proves $(i)$. To prove $(ii),$ it is sufficient to observe that $\ln \Big(\frac{r_{X}(t)}{r_{Y}(t)}\Big)$ is increasing in $t>0$ implies $\int_{0}^{t}\ln \frac{r_{X}(t)}{r_{Y}(t)}du \geq \int_{0}^{t}\ln \frac{r_{X}(u)}{r_{Y}(u)}du$ for $u\leq t$, giving  $\ln \Big(\frac{r_{X}(t)}{r_{Y}(t)}\Big)\geq \frac{1}{t}\int_{0}^{t}\ln \Big(\frac{r_{X}(u)}{r_{Y}(u)}\Big)du$ which is an equivalent condition of $X\leq_{GAI}Y$ as given in $(i)$.  Thus, $(ii)$ is proved. To prove $(iii)$, we note that $X\leq_{AF}Y$  if and only if $\ln \Big(\frac{r_{X}(t)}{r_{Y}(t)}\Big)$ is increasing in $t>0$. This implies $\ln \Big(\frac{r_{X}(t)}{r_{Y}(t)}\Big)\geq \frac{1}{t}\int_{0}^{t}\ln \Big(\frac{r_{X}(u)}{r_{Y}(u)}\Big)du,$ and this is equivalent to $X\leq_{GAI}Y$ as proved in $(i)$. This completes the proof. $\hfill\Box$\\

\hspace*{1cm}It is easy to see from the following theorem
that $FR$ order is stronger than AFR order, GFR order and HFR
order.
\begin{t1}
    If $X\leq_{FR} Y$ then $X\leq_{AFR} Y,$ $X\leq_{GFR} Y,$ and $X\leq_{HFR} Y.$
\end{t1}
{\bf Proof.} If $X\leq_{FR} Y$ then $r_{X}(t)\geq r_{Y}(t)$ for $t>0,$ which implies $\int_{0}^{t}r_{X}(u)du \geq \int_{0}^{t}r_{Y}(u)du, $ giving $X\leq_{AFR} Y.$ Similarly,  $r_{X}(t)\geq r_{Y}(t)$ for $t>0,$ implies $\int_{0}^{t}\ln r_{X}(u)du \geq \int_{0}^{t}\ln r_{Y}(u)du, $ giving $\exp\left(\frac{1}{t}\int_{0}^{t}\ln r_{X}(u)du\right)\geq \exp\left(\frac{1}{t}\int_{0}^{t}\ln r_{Y}(u)du\right).$ Thus, $X\leq_{FR} Y$ implies $X\leq_{GFR} Y.$ Similarly, $X\leq_{FR} Y$ implies $X\leq_{HFR} Y.$$\hfill\Box$\\

\hspace*{1cm}The following theorem is interesting as we note that AFR order
is equivalent to the usual stochastic order.
\begin{t1}
\label{imppp}
For two non-negative random variables $X$ and $Y,$  $X\leq_{AFR}Y$ if and only if $X\leq_{st}Y.$
\end{t1}
{\bf Proof.} Since $X\leq_{AFR}Y$ is equivalent to the fact that for all $x>0,$ $\int_{0}^{x}r_{X}(u)du\geq \int_{0}^{x}r_{Y}(u)du,$ or, equivalently, $\exp \Big(-\int_{0}^{x}r_{X}(u)du\Big)\leq \exp \Big(-\int_{0}^{x}r_{Y}(u)du\Big),$  giving $\bar{F}_{X}(t)\leq \bar{F}_{X}(t).$ This completes the proof.~$\hfill\Box$
\section{Applications in real-life data and a simulated data}
Our research extends to practical applications through the examination of a case study and  analysis of a simulated data, having relevance in the fields of survival and reliability analysis.
In this section, we take up a case study followed by simulation of a Weibull distribution.\\

We first use muhaz package (on the time points of a given data) available in Comprehensive R Archive Network to estimate $HR$. In particular, we took up epanechnikov kernel while applying muhaz. After getting estimates of $HR$ at arbitrary grid points specified by muhaz, we compute $r(\cdot), A(\cdot), G(\cdot), H(\cdot)$ and $L(\cdot),L^{G}(\cdot)$ and $L^{H}(\cdot)$ at different grid points using  RStudio 2023.09.0+463. \\

Considering the fact that $FR, AFR, GFR$ and $HFR$ are of same dimension (per unit time), they can be compared by their point-wise values. The dimension of FR is per unit time whereas $AI, GAI$ and $HAI$ are dimensionless. Despite differences in dimension their pattern can be explored for further analyses.
\subsection{Real-life data: Accelerated Life testing Data}
In the following example we explore the notion of aging function a case study of life-testing data, aiming to glean deeper insights.
\begin{ex}
\label{ex1}
Manufacturers A,B,C,D,E of hip – joint products which are made of different material combination of ball and cup of the joints were subjected to fatigue test (axial and torsional) on samples of 30 assemblies each to measure the amount of wear out particles over time. The censoring time is taken as 300 hours. The number of cycles to failure of hip joint products by each manufacturer (scaled by $10^{6}$) in accelerated life testing problem (ALT) for each test unit were noted from Elsayed (2021), Problem 5.5 page 380. A sample size of 40 units was selected out of which 20 units failed for each manufacturer.
\end{ex}
\begin{figure}
\begin{minipage}{0.5\linewidth}
		\includegraphics[width=8cm,keepaspectratio]{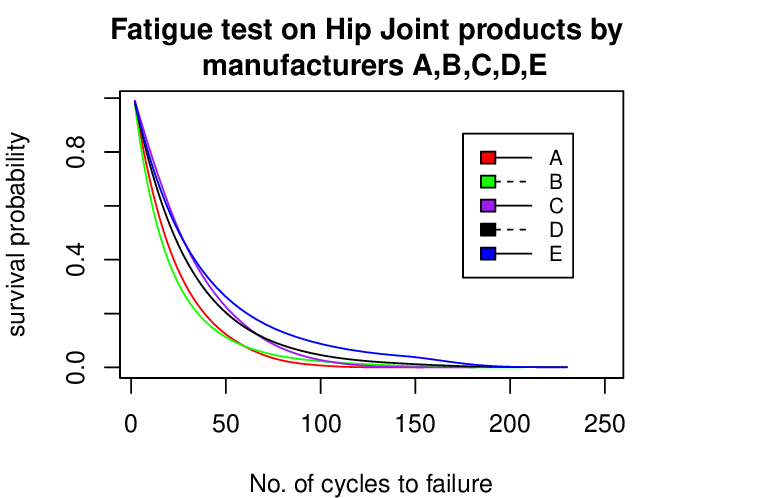}
			\end{minipage}
\hfill
		\begin{minipage}{0.5\linewidth}
		\includegraphics[width=8cm,keepaspectratio]{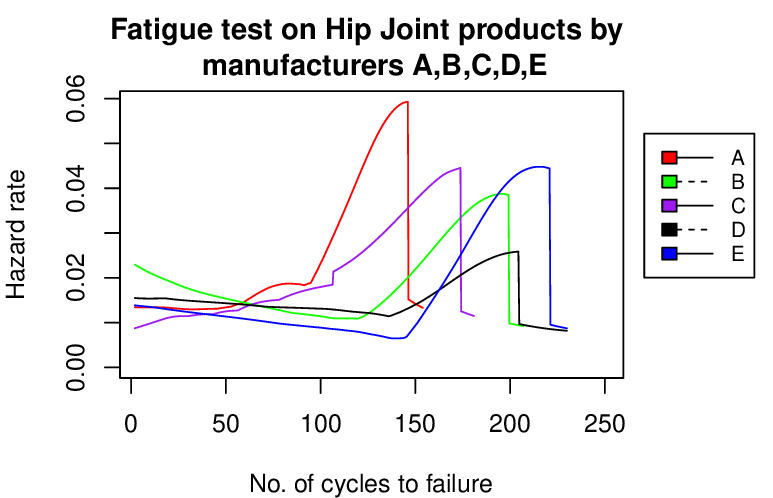}
		\end{minipage}
\hfill
	\begin{minipage}{0.5\linewidth}
		\includegraphics[width=8cm,keepaspectratio]{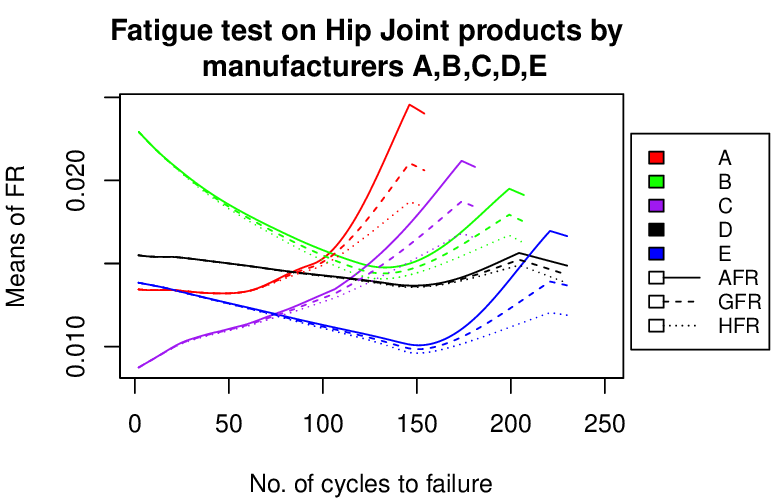}
			\end{minipage}
\begin{minipage}{0.5\linewidth}
		\includegraphics[width=8cm,keepaspectratio]{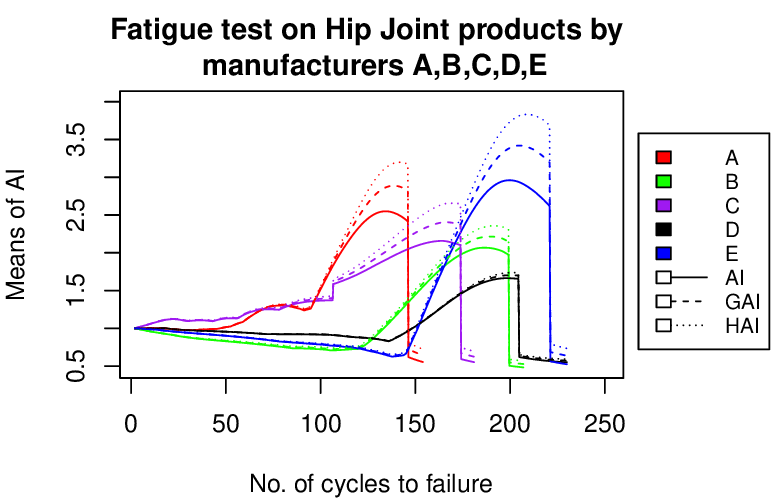}
			\end{minipage}
\\
				{{\scriptsize {\bf Figure 1}: Analysis of ALT problem through means of FR and AI (Example  \ref{ex1}).
}}
\end{figure}
\hspace*{1cm}From Figure 1, we observe that initially approximately as total number of 70 cycles to failure were observed, components from manufacturer B exhibits the highest level of efficiency, followed by manufacturer E, D, A and notably the least efficiency being attributed to component from manufacturer C.
After completion of around 110 cycles to failure, a marked decline in efficiency was observed in the components procured from manufacturer C.
Although the initial efficiency of components sourced from manufacturer D was not on par with those from manufacturer B and E, over longer period of time it demonstrated superior performance compared to all the other manufacturers.
At the outset, the components sourced from manufacturer E outperformed majority of their counterparts. However, at approximately 145th cycle to failure a conspicuous deterioration in efficiency becomes apparent, yet they do endure longer among all manufacturers.
Components from manufacturer A depict the lowest number of cycles to failure signifying depletion in their operational capacity sooner than others.\\
\hspace*{1cm}Components from Manufacturer B exhibits the highest average failure rates, namely AFR, GFR and HFR, surpassing all others, a distinction sustained until the threshold of approximately 100 cycles to failures.
The components sourced from manufacturers B, D and E exhibit an initial decline in their failure rates, followed by a gradual escalation in failure rates, indicating a diminishing level of efficiency over time.
There exists a consistent propensity for an ascending failure rate in the above graph for components acquired from manufactures A and C, suggesting continuous decrement in performance as the number of cycles to failure increases.
Despite the fact that components from manufacturer E do not commence with the most minimal failure rate initially, they progressively exhibit the least failure rate among all, establishing their status as slightly more resilient than the rest in the longer run.
At around 150th cycles to failure  components from manufacturer A depict most pronounced surge in failure rate followed by earliest diminished efficacy.\\

\hspace*{1cm}From the preceding graphical illustrations (Figure 1), we can infer that in the fatigue test products from manufacturer E can endure more numbers of cycles to failure than rest, a fact reflected through GFR and HFR in a more precise manner.
\subsection{Simulated Data}

\begin{figure}
\begin{minipage}{0.5\linewidth}
		\includegraphics[width=7cm,keepaspectratio]{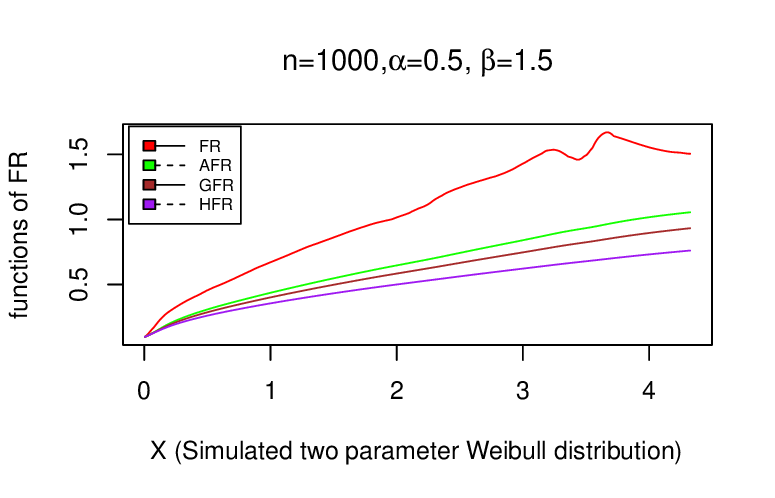}
			\end{minipage}
		\begin{minipage}{0.5\linewidth}
		\includegraphics[width=7cm,keepaspectratio]{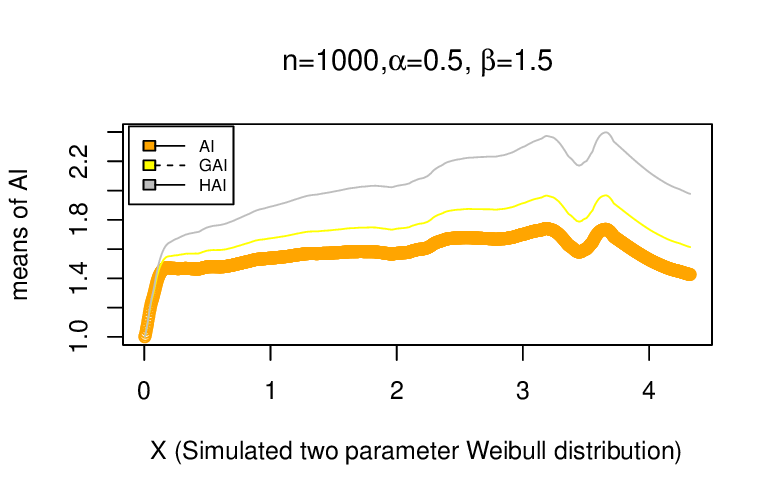}
		\end{minipage}
\hfill
	\begin{minipage}{0.5\linewidth}
		\includegraphics[width=7cm,keepaspectratio]{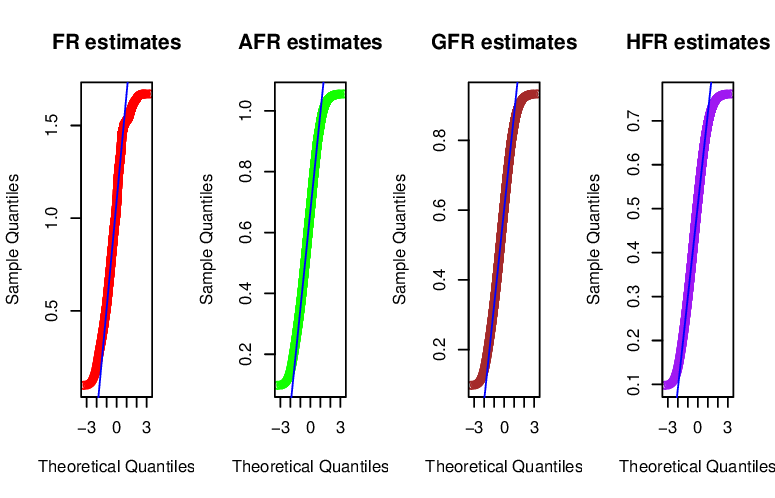}
			\end{minipage}
		\hfill
	\begin{minipage}{0.5\linewidth}
		\includegraphics[width=7cm,keepaspectratio]{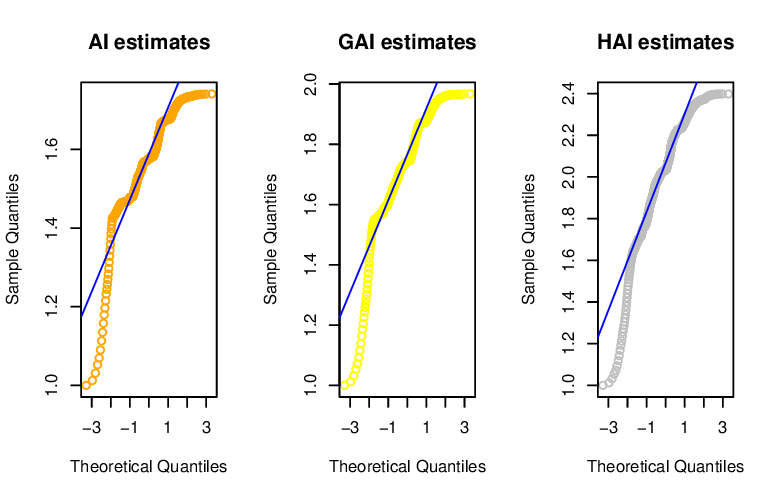}
			\end{minipage}
		\scriptsize{{\bf Figure 2:} Plot of Simulated Weibull vs. aging functions; Q-Q plot. We take $n=1000,$ sample size of Simulated \\Weibull, $\alpha=0.5,$ scale parameter; $\beta=1.5,$ shape parameter;}\\\\
\hfill
	\begin{minipage}{0.5\linewidth}
		\includegraphics[width=7cm,keepaspectratio]{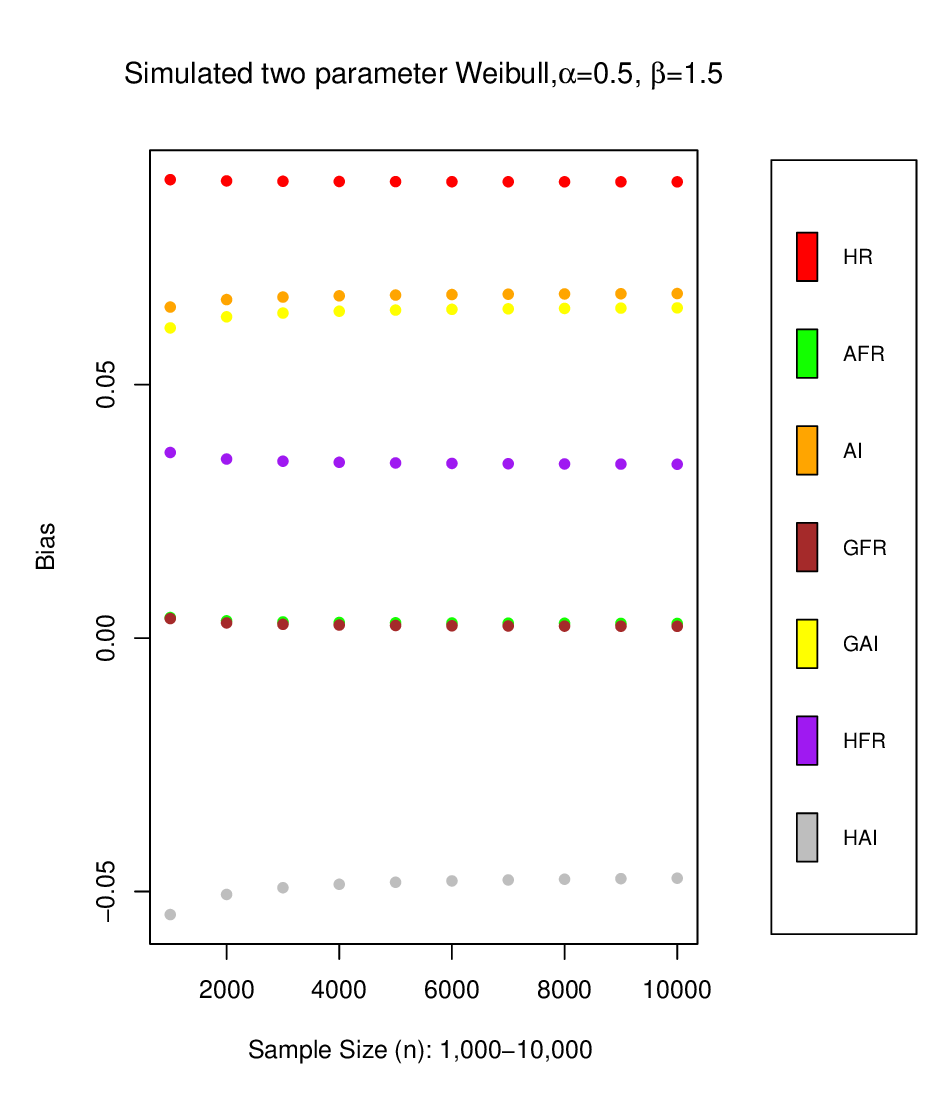}
			\end{minipage}
		\hfill
	\begin{minipage}{0.5\linewidth}
		\includegraphics[width=7cm,keepaspectratio]{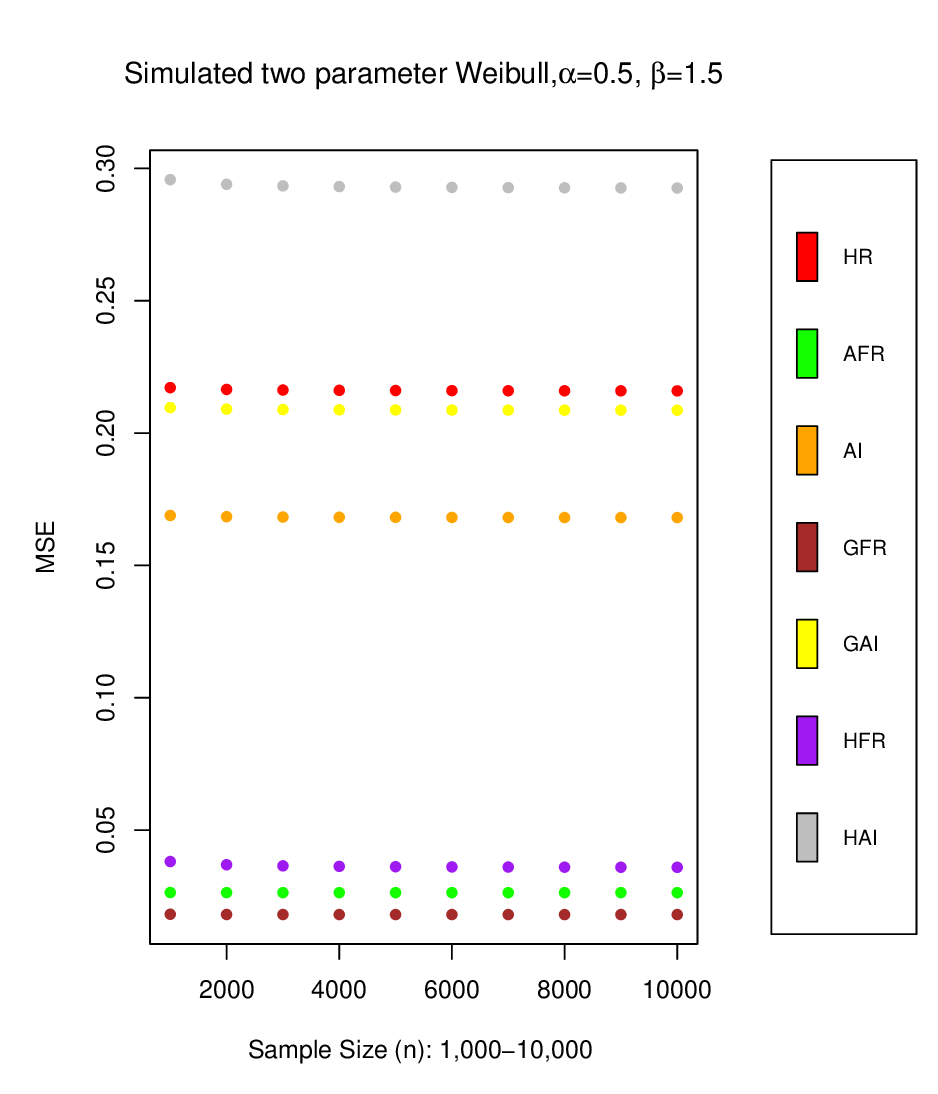}
			\end{minipage}
		{\scriptsize {\bf Figure 3:} Bias and MSE of Simulated Weibull ($\alpha=0.5$, $\beta=1.5$) with sample size ranging from 1,000 to 10,000.}
\end{figure}
\hspace*{1cm}Now, we proceed to take up a simulated study to observe aging phenomenon employing the proposed metrics as mentioned in this paper earlier. Simultaneously, in this process we also make an attempt to ascertain the measure which performs better by evaluating the bias and MSE of their corresponding estimates using simulated data.\\

\hspace*{0.5cm}We generate a sample of size $n=1000$ from two-parameter Weibull distribution having survival function $\overline{F}(t)=\exp(-\alpha t^{\beta}),$ for $t \geq 0$ with scale parameter $\alpha=0.5,$ shape parameter $\beta=1.5.$
First, we begin with obtaining the estimates of failure rate at simulated data points, subsequently that of AFR, AI, GFR, GAI, HFR, and HAI which are depicted in Figure 2 along with the respective Q-Q plot.\\

\hspace*{1cm}The actual value and estimated values of $AFR, AI, GFR, GAI, HFR,$ and $HAI$ at the simulated data points obtained from the Weibull distribution are recorded for sample sizes ranging from 1,000 to 10,000 to obtain corresponding  bias and mean square error (MSE) as illustrated in Figure 3.\\

We also infer that the bias of the estimators of aging functions  can be ordered as $HR>AI>GAI>HAI>HFR>AFR \approx GFR$. All the estimates of the aging functions have positive bias except that of $HAI$. The MSE of the corresponding estimates are found to be in the following order $HAI>HR>GAI>AI>HFR>AFR>GFR$. We note that bias of estimates of failure rate are considerably more compared to $GFR$ and $HFR.$ In fact, bias of $HR$ is larger than that of $HFR$. Further, bias of $HFR$ is greater than that of  $AFR$, whereas there is insignificant difference in bias of $AFR$ and $GFR$. A similar kind of ordering is observed among MSE of the estimates of $AFR, GFR$ and $HFR.$ \\

${}$\hspace{1cm}So one might opt for GFR as the preferred choice over AFR and HFR when analyzing the aging phenomenon, considering the fact that GFR  demonstrates least bias as well as least MSE. On the other hand, HAI exhibits the least bias but relatively higher mean squared error (MSE). In light of this, GAI may be the preferred choice over AI and HAI for the analysis of aging, given that GAI strikes a balance with moderate bias and MSE. The said inference can be applied in situations where aging phenomena can be described by two-parameter Weibull distribution.
\section{Conclusion}
${}$\hspace{1cm} In this article, we develop two new aging functions namely, Geometric Aging Intensity (GAI) and Harmonic Aging Intensity (HAI) as functions of GFR and HFR, the properties of which we explore extensively. Furthermore, we also establish relationships between the different aging intensities, namely the GAI, HAI and the existing (Arithmetic) Aging Intensity. Characterization results are presented that provide insights into specificity of the behavior of the aging functions. We also investigate their properties in connection with system reliability, specifically in the context of multi-component systems that are connected in series.\\

${}$\hspace{1cm}Extending our findings to applications with simulated and real data demonstrate how the deductive analyses derived from the graphical representation harmonize seamlessly with the theoretically proven results. The superior performance of GFR and GAI in simulations establishes them as viable alternatives to the existing measures to describe the aging phenomenon.\\

${}$\hspace{1cm}An avenue that has not been explored in this article relates to the inferential properties for the estimators of the measures. The finite and large sample properties of the estimators, both in the context of complete and censored data extend the premise of application from a limited scoped reliability setting to a more practical context.
Although the discussion in this article is confined to a simple engineering setting, the application of the measures can be explored in a biomedical context. In such a context, the study of association of such intensity functions with important covariates may be of interest. Estimation and testing of the estimators is an important future objective that is currently under investigation by the study team.\\

\section*{Disclosure statement}
No potential conflict of interest was reported by the authors.\\

\section*{Ackowledgement}
Subarna Bhattacharjee would like to thank Odisha State Higher Education Council for providing support to carry out the research project under OURIIP, Odisha, India (Grant No. 22-SF-MT-073).

\end{document}